\documentclass[a4paper,12pt]{article}
\usepackage[T1]{fontenc}
\usepackage[latin1]{inputenc}
\usepackage[english]{babel}
\usepackage{hyperref}
\usepackage[dvips]{graphics}
\usepackage{amsmath}
\usepackage{amssymb}
\usepackage{amsopn}
\usepackage{amsbsy}
\usepackage{amstext}
\usepackage{latexsym}
\usepackage{amsfonts}
\usepackage{array}
\usepackage{epsfig}
\usepackage{mathrsfs}
\usepackage{color}
\title{On the series of the reciprocals lcm's of sequences of positive integers:\\ A curious interpretation}
\author{\sc Bakir FARHI}
\date{}

\newtheorem{thm}{Theorem}
\newtheorem{prop}[thm]{Proposition}
\newtheorem{lemma}[thm]{Lemma}

\newtheorem{coll}[thm]{Corollary}

\newtheorem{thmn}{Theorem$\!\!\!$}   

\let\epsilon=\varepsilon
\def\d{{\underline{\bf{d}}}}
\def\A{{\mathscr{A}}}
\def\B{{\mathscr{B}}}
\def\F{{\mathscr{F}}}
\def\U{{\mathscr{U}}}
\def\lcm{{\rm lcm}}
\def\c{{\rm c}}
\def\EMts{\mspace{.3mu}}
\def\nb#1{{\left\vert{\EMts\EMts #1 \EMts\EMts}\right\vert}}

\def\EMdash{\leavevmode\hbox to 7.5mm{\vrule height .63ex depth -.59ex
    width 5.4mm\hfill}}

\begin{document}
\maketitle ~\vspace{-1.5cm}
\begin{center}
{\tt bakir.farhi@gmail.com}
\end{center}
\maketitle \vspace{-7cm}
\begin{flushleft}
{\it Integers: Electronic Journal of Combinatorial Number Theory}, \\
{\bf 9} (2009), p. 555-567 ($\#$A42).
\end{flushleft}~\vspace{5cm}

\begin{abstract}
In this paper, we prove the following result:
\begin{quote}
Let $\A$ be an infinite set of positive integers. For all positive
integer $n$, let $\tau_n$ denote the smallest element of $\A$
which does not divide $n$. Then we have
$$\lim_{N \rightarrow + \infty} \frac{1}{N} \sum_{n = 1}^{N} \tau_n = \sum_{n = 0}^{\infty} \frac{1}{\lcm\{a \in \A ~|~ a \leq n\}} .$$
\end{quote}
In the two particular cases when $\A$ is the set of all positive
integers and when $\A$ is the set of the prime numbers, we give a
more precise result for the average asymptotic behavior of
${(\tau_n)}_n$. Furthermore, we discuss the irrationality of the
limit of $\tau_n$ (in the average sense) by applying a result of
Erd\H{o}s.
\end{abstract}

\noindent{\bf Keywords.} Least common multiple; Special sequences
of integers; Convergence in the average sense; Irrational
numbers.~\vspace{2mm}

\noindent{\bf AMS classification.} 11B83, 40A05.

\section{Introduction and Results}

In Number Theory, it is frequent that a sequence of positive
integers does not have a regular asymptotic behavior but has a
simple and regular asymptotic average behavior. As examples, we
can cite the following:
\begin{enumerate}
\item[(i)] The sequence ${(d_n)}_{n \geq 1}$, where $d_n$ denotes
the number of divisors of $n$.
 \item[(ii)] The sequence ${(\sigma(n))}_{n \geq 1}$, where $\sigma(n)$ denotes the sum of divisors of $n$.
 \item[(iii)] The Euler totient function ${(\varphi(n))}_{n \geq 1}$, where
$\varphi(n)$ denotes the number of positive integers, not
exceeding $n$, that are relatively prime to $n$.
\end{enumerate}
We refer the reader to \cite{hw} for many other examples.\\
In this paper, we give another type of sequence which we
describe as follows:\\
Let $a_1 < a_2 < \cdots$ be an increasing sequence of positive
integers which we denote by $\A$. For all positive integers $n$,
let $\tau_n$ denote the smallest element of $\A$ which doesn't
divide $n$. Then, we shall prove the following
\begin{thm}\label{t1}
We have
$$\lim_{N \rightarrow + \infty} \frac{1}{N} \sum_{n = 1}^{N} \tau_n = \sum_{n = 0}^{\infty} \frac{1}{\lcm\{a \in \A ~|~ a \leq n\}}$$
in both cases when the series on the right-hand side converges or
diverges.
\end{thm}
In the particular cases when $\A$ is the sequence of all positive
integers and when it is the sequence of the prime numbers, we
refine the proof of Theorem \ref{t1} to obtain the following more
precise results:
\begin{coll}\label{c1}
For all positive integers $n$, let $e_n$ denote the smallest
positive integer which doesn't divide $n$. Then, we have
$$\frac{1}{N} \sum_{n = 1}^{N} e_n = \ell_1 + O_N\left(\frac{(\log N)^2}{N \log\log N}\right) ,$$
where
$$\ell_1 := \sum_{n \in \mathbb{N}}\frac{1}{\lcm(1 , 2 , \dots , n)} < + \infty .$$
\end{coll}
\begin{coll}\label{c2}
For all positive integer $n$, let $q_n$ denote the smallest prime
number which doesn't divide $n$. Then, we have
$$\frac{1}{N} \sum_{n = 1}^{N} q_n = \ell_2 + O_N\left(\frac{(\log N)^2}{N \log\log N}\right) ,$$
where
$$\ell_2 := \sum_{n \in \mathbb{N}}\frac{1}{~~~~\prod_{{}_{{}_{\!\!\!\!\!\!\!\!\!\!\!\!\!\!\!\!\!\!\!\!\text{$p$ prime, $p \leq n$}}}} p} < + \infty .$$
\end{coll}
Further, by applying a result of Erd\H{o}s \cite{erd}, we derive a
sufficient condition for the average limit of ${(\tau_n)}_n$ to be
an irrational number. We have the following.
\begin{prop}\label{p1}
Let $\d(\A)$ denote the lower asymptotic density of $\A$, that is
$$\d(\A) := \liminf_{N \rightarrow + \infty} \frac{1}{N}
\sum_{\begin{subarray}{c} \scriptstyle{a \in \A} \\ \scriptstyle{a
\leq N}
\end{subarray}} 1 .$$
Suppose that $\d(\A) > 1 - \log{2}$. Then, the average limit of
${(\tau_n)}_n$ is an irrational number.\\
In particular, the numbers $\ell_1$ and $\ell_2$ appearing
respectively in Corollaries \ref{c1} and \ref{c2} are irrational.
\end{prop}
\section{The Proofs}
\subsection{Some preparations and preliminary results}

Throughout this paper, we let $\mathbb{N}^*$ denote the set
$\mathbb{N} \setminus \{0\}$ of all positive integers. For a given
real number $x$, we let $\lfloor x \rfloor$ and $\langle x
\rangle$ denote respectively the integer part and the fractional
part of $x$. Further, we adopt the natural convention that the
least common multiple and the product of the elements of an empty
set are equal to $1$.

We fixe an increasing sequence of positive integers $a_1 < a_2 <
\cdots$, which we denote by $\A$ and for all positive integer $n$,
we let $\tau_n$ denote the smallest element of $\A$ which doesn't
divide $n$.

For all $\alpha \in \A$, we let $L(\alpha)$ denote the positive
integer defined by
\begin{equation}\label{eq1}
L(\alpha) := \frac{\lcm\{a \in \A ~|~ a \leq \alpha\}}{\lcm\{a \in
\A ~|~ a < \alpha\}} .
\end{equation}
We let then $\B$ denote the subset of $\A$ defined by
\begin{equation}\label{eq2}
\B := \left\{a \in \A ~|~ L(a) > 1\right\} .
\end{equation}
We shall see later that $\B$ is just the set of the values of the
sequence ${(\tau_n)}_n$. We begin with the following lemma
\begin{lemma}\label{l1}
For all positive integer $n$, we have
$$\lcm\left\{a \in \A ~|~ a \leq n\right\} = \lcm\left\{b \in \B ~|~ b \leq n\right\} .$$
\end{lemma}

\noindent{\bf Proof.} Let $n \geq 1$ and let $a_1 , \dots , a_k$
be the elements of $\A$ not exceeding $n$. By using the following
well-known property of the least common multiple:
$$\lcm(a_1 , \dots , a_i , \dots , a_k) = \lcm(\lcm(a_1 , \dots , a_i) , a_{i + 1} , \dots , a_k) ~~~~ \text{(for $i = 1 , 2 , \dots , k$)} ,$$
we remark that when $a_i \not\in \B$, we have $L(a_i) = 1$ and
then $\lcm(a_1 , \dots , a_i) = \lcm(a_1 , \dots , a_{i - 1})$.
So, each $a_i$ not belonging to $\B$ can be eliminated from the
list $a_1 , \dots , a_k$ without changing the value of the least
common multiple of that list. The Lemma follows.\hfill
$\blacksquare$~\vspace{2mm}

From Lemma \ref{l1}, we derive another formula for $L(\alpha)$
$(\alpha \in \A)$. We have for all $\alpha \in \A$
\begin{equation}\label{eq19}
L(\alpha) = \frac{\lcm\{b \in \B ~|~ b \leq \alpha\}}{\lcm\{b \in
\B ~|~ b < \alpha\}} .
\end{equation}

The next lemma gives a useful characterization for the terms of
the sequence ${(\tau_n)}_n$.
\begin{lemma}\label{l2}
For all positive integers $n$ and all $\alpha \in \A$, we have:
$$\tau_n = \alpha \Longleftrightarrow \exists k \in \mathbb{N}^* , k \not\equiv 0 \!\!\!\!\mod L(\alpha)
~\text{such that}~ n = k \cdot \lcm\{b \in \B ~|~ b < \alpha\} .$$
\end{lemma}

\noindent{\bf Proof.} Let $n \in \mathbb{N}^*$ and $\alpha \in
\A$. By the definition of the sequence ${(\tau_n)}_n$, the
equality $\tau_n = \alpha$ amounts to saying that $n$ is a
multiple of each element $a \in \A$ satisfying $a < \alpha$ and
that $n$ is not a multiple of $\alpha$. Equivalently, $\tau_n =
\alpha$ if and only if $n$ is a multiple of $\lcm\{a \in \A ~|~ a
< \alpha\}$ without being a multiple of $\lcm\{a \in \A ~|~ a \leq
\alpha\}$. So, it suffices to set $k := \frac{n}{\lcm\{a \in \A
~|~ a < \alpha\}}$ to obtain the equivalence:
$$\tau_n = \alpha \Longleftrightarrow \exists k \in \mathbb{N}^* , k \not\equiv 0 \!\!\!\!\mod L(\alpha) ~\text{such that}~
n = k \cdot \lcm\{a \in \A ~|~ a < \alpha\} .$$ The lemma then
follows from Lemma \ref{l1}.\hfill $\blacksquare$

\begin{coll}\label{r1} The sequence ${(\tau_n)}_n$ takes
its values in the set $\B$. Besides, any element of $\B$ is taken
by ${(\tau_n)}_n$ infinitely often.
\end{coll}

\noindent{\bf Proof.} Let $\alpha \in \A$. If $\alpha \not\in \B$,
then we have $L(\alpha) = 1$ and thus there is no $k \in
\mathbb{N}^*$ such that $k \not\equiv 0 \mod L(\alpha)$. It
follows, according to Lemma \ref{l2}, that $\alpha$ cannot be a
value of ${(\tau_n)}_n$.\\
Next, if $\alpha \in \B$, then $L(\alpha) \geq 2$ and thus there
are infinitely many $k \in \mathbb{N}^*$ such that $k \not\equiv 0
\mod L(\alpha)$. This implies (according to Lemma \ref{l2}) that
there are infinitely many $n \in \mathbb{N}^*$ satisfying $\tau_n
= \alpha$. The corollary is proved.\hfill
$\blacksquare$~\vspace{2mm}

Actually, given $\alpha \in \B$, Lemma \ref{l2} even gives an
estimation for the number of solutions of the equation $\tau_n =
\alpha$ in an interval $[1 , x]$ $(x \in \mathbb{R}^+)$. For $x >
0$ and $\alpha \in \B$, define
$$\varphi(\alpha ; x) := \#\{n \in \mathbb{N}^* , n \leq x ~|~ \tau_n = \alpha\} .$$
Then, we have the following:
\begin{coll}\label{c3}
Let $\alpha \in \B$ and $x > 0$. Then we have
$$\varphi(\alpha ; x) = \frac{L(\alpha) - 1}{\lcm\{b \in \B ~|~ b \leq \alpha\}} . x + \c_{\alpha , x} ,$$
where $|\c_{\alpha , x}| < 1$. Furthermore, if $\lcm\{b \in \B ~|~
b < \alpha\} > x$, then $\varphi(\alpha ; x) = 0$.
\end{coll}

\noindent{\bf Proof.} Let $\alpha \in \B$ and $x > 0$. By Lemma
\ref{l2}, we have
\begin{eqnarray*}
\varphi(\alpha ; x) & = & \#\left\{k \in \mathbb{N}^* ~|~ k \leq
\frac{x}{\lcm\{b \in \B ~|~ b < \alpha\}} ~~\text{and}~~ k
\not\equiv 0 \!\!\!\!\mod L(\alpha)\right\} \\
& = & \left\lfloor\frac{x}{\lcm\{b \in \B ~|~ b <
\alpha\}}\right\rfloor - \left\lfloor\frac{x}{L(\alpha) . \lcm\{b
\in \B ~|~ b < \alpha\}}\right\rfloor \\
& = & \frac{x}{\lcm\{b \in \B ~|~ b < \alpha\}} -
\frac{x}{L(\alpha) . \lcm\{b \in \B ~|~ b < \alpha\}} + \c_{\alpha
, x} ,
\end{eqnarray*}
where $\c_{\alpha , x} := - \langle \frac{x}{\lcm\{b \in \B ~|~ b
< \alpha\}}\rangle + \langle\frac{x}{L(\alpha) . \lcm\{b \in \B
~|~ b < \alpha\}}\rangle$. So, it is clear that $|\c_{\alpha , x}|
< 1$. Next, we have
\begin{equation*}
\begin{split}
\frac{x}{\lcm\{b \in \B ~|~ b < \alpha\}} - \frac{x}{L(\alpha) .
\lcm\{b \in \B ~|~ b < \alpha\}} &= \frac{L(\alpha) -
1}{L(\alpha) . \lcm\{b \in \B ~|~ b < \alpha\}} . x \\
&= \frac{L(\alpha) - 1}{\lcm\{b \in \B ~|~ b \leq \alpha\}} . x ,
\end{split}
\end{equation*}
by using (\ref{eq19}). This confirms the first part of the
corollary.

Now, if $\lcm\{b \in \B ~|~ b < \alpha\} > x$, then none of the
integers of the range $[1 , x]$ is a multiple of $\lcm\{b \in \B
~|~ b < \alpha\}$. It follows, according to Lemma \ref{l2}, that
the equation $\tau_n = \alpha$ doesn't have any solution in the
range $[1 , x]$. Hence $\varphi(\alpha , x) = 0$.\\
This confirms the second part of the corollary and ends this
proof.\hfill $\blacksquare$~\vspace{2mm}

Now, let $b_1 < b_2 < \cdots$ be the elements of $\B$. To prove
our main result, we will need some properties of the sequence
${(b_n)}_n$. For simplicity, let for all $n \geq 1$
$$L_n := L(b_n) \geq 2 .$$
Then, by (\ref{eq19}), we have
$$L_n = \frac{\lcm\{b \in \B ~|~ b \leq b_n\}}{\lcm\{b \in \B ~|~ b <
b_n\}} = \frac{\lcm(b_1 , b_2 , \dots , b_n)}{\lcm(b_1 , b_2 ,
\dots , b_{n - 1})} ,$$ which gives
$$\lcm(b_1 , b_2 , \dots , b_n) = L_n \cdot \lcm(b_1 , b_2 , \dots , b_{n - 1}) ~~~~~~ (\forall n \geq 1) .$$
By iteration, we obtain for all $n \geq 1$
\begin{equation}\label{eq4}
\lcm(b_1 , b_2 , \dots , b_n) = L_1 L_2 \cdots L_n \geq 2^n .
\end{equation}

For the simplicity of some formulas in what follows, it is useful
to set $b_0 := 0$ and $L_0 := 1$. Note that $b_0 \not\in \B$. We
have the two following lemmas:
\begin{lemma}\label{l3}
We have
$$\sum_{n \in \mathbb{N}}\frac{1}{\lcm\{a \in \A ~|~ a \leq n\}} = \sum_{k = 1}^{\infty}\frac{b_k - b_{k - 1}}{L_1 L_2 \cdots L_{k - 1}} .$$
\end{lemma}

\noindent{\bf Proof.} According to Lemma \ref{l1}, we have
\begin{eqnarray*}
\sum_{n \in \mathbb{N}}\frac{1}{\lcm\{a \in \A ~|~ a \leq n\}} & =
& \sum_{n \in \mathbb{N}}\frac{1}{\lcm\{b \in \B ~|~ b \leq n\}} \\
& = & \sum_{k = 1}^{\infty}\sum_{b_{k - 1} \leq n <
b_k}\frac{1}{\lcm\{b \in \B ~|~ b \leq n\}} \\
& = & \sum_{k = 1}^{\infty}\frac{b_k - b_{k - 1}}{\lcm(b_1 , b_2
, \dots , b_{k - 1})} \\
& = & \sum_{k = 1}^{\infty}\frac{b_k - b_{k - 1}}{L_1 L_2 \cdots
L_{k - 1}} ~~~~~~~~~~~~\text{(according to (\ref{eq4})).}
\end{eqnarray*}
The Lemma is proved.\hfill $\blacksquare$

\begin{lemma}\label{l4}
Let $r$ be a positive integer and $N$ be an integer such that
$$\lcm(b_1 , b_2 , \dots , b_{r - 1}) \leq N < \lcm(b_1 , b_2 , \dots , b_r) .$$
Then, we have
$$\frac{1}{N}\sum_{n = 1}^{N} \tau_n = S_1(N) + S_2(N) f(N , r) ,$$
where
\begin{equation*}
S_1(N) := - \frac{b_r}{L_1 L_2 \cdots L_r} + \sum_{k =
1}^{r}\frac{b_k - b_{k - 1}}{L_1 L_2 \cdots L_{k - 1}} ~~,~~
S_2(N) := \frac{1}{N}\sum_{k = 1}^{r} b_k
\end{equation*}
and $f(N , r)$ is a function of $N$ and $r$, satisfying $|f(N ,
r)| < 1$.
\end{lemma}

\noindent{\bf Proof.} According to Corollaries \ref{r1}, \ref{c3}
and to the relation (\ref{eq4}), we have
\begin{eqnarray*}
\frac{1}{N}\sum_{n = 1}^{N} \tau_n & = & \frac{1}{N}\sum_{\alpha
\in \B} ~ \sum_{1 \leq n \leq N , \tau_n = \alpha} \alpha \\
& = & \frac{1}{N}\sum_{\alpha \in \B} \alpha \varphi(\alpha ; N)
\\
& = & \frac{1}{N}\sum_{\begin{subarray}{c} \scriptstyle{\alpha \in
\B}
\\ \scriptstyle{\lcm\{b \in \B ~|~ b < \alpha\} \leq N}
\end{subarray}} \!\!\!\!\!\!\!\!\!\!\!\!\!\!\!\!\alpha \varphi(\alpha ; N) \\
& = & \frac{1}{N}\sum_{k = 1}^{r} b_k \varphi(b_k ; N) \\
& = & \frac{1}{N} \sum_{k = 1}^{r} b_k \left(\frac{L_k - 1}{L_1
L_2 \cdots L_k} N + \c_{k , N}\right) ~~~~~~~~~~~~~~ \text{(where $\nb{\c_{k ,N}} < 1$)}\\
& = & \sum_{k = 1}^{r} b_k \left(\frac{1}{L_1 L_2 \cdots L_{k -
1}} - \frac{1}{L_1 L_2 \cdots L_k}\right) + \frac{1}{N} \sum_{k =
1}^{r} b_k \c_{k , N} .
\end{eqnarray*}
Then, the lemma follows by remarking that
\begin{eqnarray*}
\sum_{k = 1}^{r} b_k & ~ &
\!\!\!\!\!\!\!\!\!\!\!\!\!\!\!\left(\frac{1}{L_1 L_2 \cdots L_{k -
1}} - \frac{1}{L_1 L_2 \cdots L_k}\right) \\
& = & \sum_{k = 1}^{r}\left(\frac{b_{k - 1}}{L_1 L_2 \cdots L_{k -
1}} - \frac{b_k}{L_1 L_2 \cdots L_k}\right) + \sum_{k = 1}^{r}
\frac{b_k - b_{k - 1}}{L_1 L_2 \cdots L_{k - 1}} \\
& = & - \frac{b_r}{L_1 L_2 \cdots L_r} + \sum_{k = 1}^{r}
\frac{b_k - b_{k - 1}}{L_1 L_2 \cdots L_{k - 1}}
\end{eqnarray*}
and by defining
$$f(N , r) := \frac{\sum_{k = 1}^{r} b_k \c_{k , N}}{\sum_{k = 1}^{r} b_k} ,$$
which satisfies $\nb{f(N , r)} < 1$ (since $\nb{\c_{k , N}} < 1$
for all $k \geq 1$). This completes the proof.\hfill
$\blacksquare$~\vspace{2mm}

We will finally need two lemmas on the convergence of sequences
and series.
\begin{lemma}\label{c4}
Let ${(x_n)}_{n \geq 1}$ be a real non-increasing sequence.
Suppose that the series $\sum_{n = 1}^{\infty} x_n$ converges.
Then, we have
$$x_n = o\left(\frac{1}{n}\right) ~~~~~~ \text{(as $n$ tends to infinity)} .$$
\end{lemma}
\begin{lemma}\label{l6}
Let ${(\theta_n)}_{n \geq 1}$ be a sequence of real numbers.
Suppose that $\theta_n$ tends to $0$ as $n$ tends to infinity.
Then the sequence with general term given by
$$\frac{1}{2^n}\sum_{k = 1}^{n} 2^k \theta_k$$
also tends to $0$ as $n$ tends to infinity.
\end{lemma}

\noindent{\bf Remark.} Note that Lemma \ref{l6} is in fact a
particular case of a more general theorem in summability theory,
called Silverman--Toeplitz theorem (see e.g. \cite{h}).

\subsection{Proofs of the main results}

\noindent{\bf Proof of Theorem \ref{t1}.} Let $N$ be a positive
integer. Since (according to (\ref{eq4})) the sequence ${(\lcm(b_1
, b_2 , \dots , b_n))}_n$ increases and tends to infinity with
$n$, $N$ must lie somewhere between two consecutive terms of this
sequence. So, let $r \geq 1$ such that
$$\lcm(b_1 , b_2 , \dots , b_{r - 1}) \leq N < \lcm(b_1 , b_2 , \dots , b_r) .$$
Then, by using Lemma \ref{l4}, we have
\begin{equation}\label{eq6}
\frac{1}{N}\sum_{n = 1}^{N} \tau_n = S_1(N) + S_2(N) f(N) ,
\end{equation}
where
\begin{eqnarray}
S_1(N) & := & - \frac{b_r}{L_1 L_2 \cdots L_r} + \sum_{k =
1}^{r}\frac{b_k - b_{k - 1}}{L_1 L_2 \cdots L_{k - 1}} , \label{N1} \\
S_2(N) & := & \frac{1}{N}\sum_{k = 1}^{r} b_k \label{N2}
\end{eqnarray}
and $|f(N)| <  1$.\\
Next, set
$$S := \sum_{n = 0}^{\infty}\frac{1}{\lcm\{a \in \A ~|~ a \leq n\}} .$$
To prove Theorem \ref{t1}, we distinguish two cases according to
whether $S$ converges or diverges.

{\bf 1\textsuperscript{st} case:} $\boldsymbol{S < + \infty}${\bf
.} In this case, since the sequence ${(1/\lcm\{a \in \A ~|~ a <
n\})}_{n \geq 1}$ is clearly non-increasing, then by Corollary
\ref{c4}, we have
$$\lim_{n \rightarrow + \infty} \frac{n}{\lcm\{a \in \A ~|~ a < n\}} = 0 .$$
By specializing in this limit $n$ to the integers $b_k$ $(k \geq
1)$, we obtain (according to Lemma \ref{l1} and formula
(\ref{eq4})) that
\begin{equation}\label{eq7}
\lim_{k \rightarrow + \infty} \frac{b_k}{L_1 L_2 \cdots L_{k - 1}}
= 0 .
\end{equation}
On the one hand, according to (\ref{eq7}) and to Lemma \ref{l3},
we have (because $r$ tends to infinity with $N$)
\begin{equation}\label{eq8}
\lim_{N \rightarrow + \infty} S_1(N) = \sum_{k = 1}^{\infty}
\frac{b_k - b_{k - 1}}{L_1 L_2 \cdots L_{k - 1}} = S
\end{equation}
and on the other hand we have
\begin{eqnarray*}
S_2(N) & := & \frac{1}{N}\sum_{k = 1}^{r} b_k \\
& \leq & \frac{1}{\lcm(b_1 , \dots , b_{r - 1})} \sum_{k = 1}^{r}
b_k \\
& = & \sum_{k = 1}^{r} \frac{b_k}{L_1 L_2 \cdots L_{r - 1}} \\
& \leq & \sum_{k = 1}^{r} \frac{b_k}{L_1 L_2 \cdots L_{k - 1}}
2^{k
- r} ~~~~ \text{(since $L_i \geq 2$ for all $i \geq 1$)} \\
& = & \frac{1}{2^r} \sum_{k = 1}^{r} 2^k \frac{b_k}{L_1 L_2 \cdots
L_{k - 1}} .
\end{eqnarray*}
But by applying Lemma \ref{l6} for $\theta_k := \frac{b_k}{L_1 L_2
\cdots L_{k - 1}}$ which is seen (from (\ref{eq7})) to tend to $0$
as $k$ tends to infinity, we have
$$\lim_{r \rightarrow + \infty}\frac{1}{2^r}\sum_{k = 1}^{r} 2^k \frac{b_k}{L_1 L_2 \cdots L_{k - 1}} = 0 .$$
So, it follows (because $r$ tends to infinity with $N$) that
\begin{equation}\label{eq9}
\lim_{N \rightarrow + \infty} S_2(N) = 0 .
\end{equation}
Finally, by inserting (\ref{eq8}) and (\ref{eq9}) into
(\ref{eq6}), we get
$$\lim_{N \rightarrow + \infty} \frac{1}{N} \sum_{n = 1}^{N} \tau_n = S ,$$
as required.

{\bf 2\textsuperscript{nd} case:} $\boldsymbol{S = + \infty}${\bf
.} In this case, by using (\ref{eq6}), we are going to bound from
below $\frac{1}{N} \sum_{n = 1}^{N} \tau_n$ by an expression
tending to infinity with $N$.\\
\indent On the one hand, we have
\begin{eqnarray}
S_1(N) & := & - \frac{b_r}{L_1 L_2 \cdots L_r} + \sum_{k = 1}^{r}
\frac{b_k - b_{k - 1}}{L_1 L_2 \cdots L_{k - 1}} \notag \\
& = & - \frac{b_r}{L_1 L_2 \cdots L_r} + \frac{b_r - b_{r -
1}}{L_1 L_2 \cdots L_{r - 1}} + \sum_{k = 1}^{r - 1} \frac{b_k -
b_{k -
1}}{L_1 L_2 \cdots L_{k - 1}} \notag \\
& \geq & \frac{b_r}{L_1 L_2 \cdots L_{r - 1}} - 2 + \sum_{k =
1}^{r - 1} \frac{b_k - b_{k - 1}}{L_1 L_2 \cdots L_{k - 1}}
\label{eq10}
\end{eqnarray}
(because we obviously have $b_i \leq \lcm(b_1 , \dots , b_i) = L_1
L_2 \cdots L_i$, for all $i \geq 1$). \\
\indent On the other hand, we have
\begin{eqnarray}
S_2(N) & := & \frac{1}{N} \sum_{k = 1}^{r} b_k \notag \\
& \leq & \sum_{k = 1}^{r} \frac{b_k}{L_1 L_2 \cdots L_{r - 1}} \notag \\
& \leq & \frac{b_r}{L_1 L_2 \cdots L_{r - 1}} + \sum_{k = 1}^{r -
1} \frac{L_1 L_2 \cdots L_k}{L_1 L_2 \cdots L_{r - 1}} ~~ \text{(since $b_i \leq L_1 L_2 \cdots L_i$, for all $i$)} \notag \\
& = & \frac{b_r}{L_1 L_2 \cdots L_{r - 1}} + \sum_{k = 1}^{r - 1}
\frac{1}{L_{k + 1} L_{k + 2} \cdots L_{r - 1}} \notag \\
& \leq & \frac{b_r}{L_1 L_2 \cdots L_{r - 1}} + \sum_{k = 1}^{r -
1} \frac{1}{2^{r - k - 1}} ~~~~~~~~~~~~ \text{(since $L_i \geq 2$
for all
$i$)} \notag \\
& < & \frac{b_r}{L_1 L_2 \cdots L_{r - 1}} + 2 . \label{eq11}
\end{eqnarray}
It follows, by inserting (\ref{eq10}) and (\ref{eq11}) into
(\ref{eq6}), that
\begin{eqnarray*}
\frac{1}{N} \sum_{n = 1}^{N} \tau_n & = & S_1(N) + S_2(N) f(N) \\
& \geq & S_1(N) - S_2(N) ~~~~ \text{(since $|f(N)| < 1$)} \\
& \geq & \sum_{k = 1}^{r - 1} \frac{b_k - b_{k - 1}}{L_1 L_2
\cdots L_{k - 1}} - 4 .
\end{eqnarray*}
But since (according to Lemma \ref{l3}) $\sum_{k = 1}^{\infty}
\frac{b_k - b_{k - 1}}{L_1 L_2 \cdots L_{k - 1}} = S = + \infty$,
we conclude that
$$\lim_{N \rightarrow + \infty} \frac{1}{N} \sum_{n = 1}^{N} \tau_n = + \infty .$$
This completes the proof of the theorem.\hfill
$\blacksquare$~\vspace{2mm}

\noindent {\bf Proof of Corollary \ref{c1}.} In the situation of
Corollary \ref{c1}, $\A$ is the set of all positive integers and
then $\B$ is the set of the powers of prime numbers. We must
repeat the proof of Theorem \ref{t1} and give more precision to
the two quantities $\frac{b_r}{L_1 L_2 \cdots L_r}$ and $\sum_{k =
1}^{r} b_k$. First let us show that
\begin{equation}\label{eq12}
b_r \sim \log N ~~~~~~ \text{(as $N$ tends to infinity).}
\end{equation}
By the definition of $r$, recall that
\begin{equation}\label{eq13}
\lcm(b_1 , b_2 , \dots , b_{r - 1}) \leq N < \lcm(b_1 , b_2 ,
\dots , b_r) .
\end{equation}
Next, by Lemma \ref{l1}, we have (since $\A = \mathbb{N}^*$)
\begin{equation*}
\begin{split}
\lcm(b_1 , b_2 , \dots , b_r) &= \lcm(1 , 2 , \dots , b_r) , \\
\lcm(b_1 , b_2 , \dots , b_{r - 1}) &= \lcm(1 , 2 , \dots , b_{r -
1})
\end{split}
\end{equation*}
and by the prime number theorem (see e.g. \cite{hw}), we have on
the one hand
$$
\log\lcm(1 , 2 , \dots , n) \sim n ~~~~~~ \text{(as $n$ tends to
infinity)}
$$
and on the other hand (because $\B$ is the set of the powers of prime
numbers)
$$
b_n \sim b_{n - 1} ~~~~~~ \text{(as $n$ tends to infinity).}
$$
Taking into account all these facts, we derive from (\ref{eq13})
that effectively $b_r \sim \log N$, confirming (\ref{eq12}).

Now, we are going to be precise the order of magnitude of
$\frac{b_r}{L_1 L_2 \cdots L_r}$ and $\sum_{k = 1}^{r} b_k$. We
have
$$\frac{b_r}{L_1 L_2 \cdots L_r} = \frac{b_r}{\lcm(b_1 , b_2 , \dots , b_r)} < \frac{b_r}{N} ~~~~~~ \text{(according to (\ref{eq13}))} .$$
It follows, according to (\ref{eq12}), that
\begin{equation}\label{eq14}
\frac{b_r}{L_1 L_2 \cdots L_r} = O\left(\frac{\log N}{N}\right) .
\end{equation}
Next, because $\B$ is the set of the powers of the prime numbers,
we have
$$
\sum_{k = 1}^{r} b_k = \sum_{\begin{subarray}{c} \scriptstyle{e
\geq 1 , p ~\text{prime}} \\ \scriptstyle{p^e \leq b_r}
\end{subarray}} \!\!\!\!\!\!\!\!p^e = \sum_{\begin{subarray}{c} \scriptstyle{p ~\text{prime}} \\ \scriptstyle{p \leq b_r}
\end{subarray}}\left(\phantom{\sum}\right.\!\!\!\!\!\!\!\!\!\sum_{1 \leq e \leq \lfloor\frac{\log b_r}{\log p}\rfloor}
p^e\!\!\!\!\!\!\!\!\!\left.\phantom{\sum}\right) =
\!\!\!\!\sum_{\begin{subarray}{c} \scriptstyle{p ~\text{prime}}
\\ \scriptstyle{p \leq b_r}
\end{subarray}}\frac{p}{p - 1}\left(p^{\lfloor\frac{\log b_r}{\log p}\rfloor} - 1\right)
.$$ Since for any prime number $p$, we have $\frac{p}{p - 1} \leq
2$ and $p^{\lfloor \log b_r / \log p\rfloor} \leq p^{\log b_r /
\log p} = b_r$, it follows that
$$\sum_{k = 1}^{r} b_k \leq 2 b_r \pi(b_r) ,$$
where $\pi$ denotes the prime-counting function. But, by using
again the prime number theorem and (\ref{eq12}), we have $b_r
\pi(b_r) = O(\frac{b_r^2}{\log b_r}) = O(\frac{(\log
N)^2}{\log\log N})$. Hence
\begin{equation}\label{eq15}
\sum_{k = 1}^{r} b_r = O\left(\frac{(\log N)^2}{\log\log N}\right)
.
\end{equation}
It finally remains to insert (\ref{eq14}) and (\ref{eq15}) into
(\ref{N1}) and (\ref{N2}) respectively to obtain (according to
(\ref{eq6}) and to Lemma \ref{l3}) that
$$
\lim_{N \rightarrow + \infty} \frac{1}{N} \sum_{n = 1}^{N} \tau_n
= \sum_{n \in \mathbb{N}} \frac{1}{\lcm(1 , 2 , \dots , n)} +
O\left(\frac{(\log N)^2}{N \log\log N}\right) .
$$
The corollary is proved.\hfill $\blacksquare$~\vspace{2mm}

\noindent {\bf Proof of Corollary \ref{c2}.} In the situation of
Corollary \ref{c2}, $\A$ is the set of the prime numbers and then
$\B = \A$. So, for all $n \geq 1$, we have $a_n = b_n = L_n =
p_n$, where $p_n$ denotes the $n$\textsuperscript{th} prime
number. Consequently, we have (in the context of the proof of
Theorem \ref{t1})
\begin{equation}\label{eq16}
p_1 p_2 \cdots p_{r - 1} \leq N < p_1 p_2 \cdots p_r .
\end{equation}
So, by the prime number theorem (see e.g. \cite{hw}), we have
$$p_r \sim \log N ~~~~~~ \text{(as $N$ tends to infinity).}$$
From this last, it follows that
$$\frac{b_r}{L_1 L_2 \cdots L_r} = \frac{p_r}{p_1 p_2 \cdots p_r} < \frac{p_r}{N} \sim \frac{\log N}{N} ~~~~~~~~~~\text{(as $N$ tends to infinity),}$$
which gives
\begin{equation}\label{eq17}
\frac{b_r}{L_1 L_2 \cdots L_r} = O\left(\frac{\log N}{N}\right)
\end{equation}
and that
$$\sum_{k = 1}^{r} b_k = \sum_{k = 1}^{r} p_k \sim \frac{p_r^2}{2 \log p_r} \sim \frac{(\log N)^2}{2 \log\log N}
~~~~~~~~~~\text{(as $N$ tends to infinity),}$$ which gives
\begin{equation}\label{eq18}
\frac{1}{N} \sum_{k = 1}^{r} b_k = O\left(\frac{(\log N)^2}{N
\log\log N}\right) .
\end{equation}
To conclude, it suffices to insert (\ref{eq17}) and (\ref{eq18})
into (\ref{N1}) and (\ref{N2}) respectively and use (\ref{eq6})
and Lemma \ref{l3}. The result of the corollary follows.\hfill
$\blacksquare$~\vspace{2mm}

Now, we are going to prove Proposition \ref{p1}. To do so, we need
the following result of Erd\H{o}s \cite{erd}.

\begin{thmn}[Erd\H{o}s \cite{erd}]
Let $u_1 < u_2 < \cdots$ be an infinite sequence of positive
integers. Set $\U := \{u_1 , u_2 , \dots\}$ and suppose that
$$\d(\U) > 1 - \log 2 = 0.306\dots$$
Then, the real positive number
$$\sum_{n = 1}^{\infty} \frac{1}{\lcm\{u \in \U ~|~ u \leq n\}}$$
is irrational.
\end{thmn}

\noindent{\bf Proof of Proposition \ref{p1}.} The first part of
Proposition \ref{p1} which concerns a general set $\A$ is clearly
an immediate consequence of the Main Theorem \ref{t1} and the
above theorem of Erd\H{o}s. Next, since the set $\mathbb{N}^*$ of
all positive integers has asymptotic density $1 > 1 - \log 2$, the
irrationality of the constant $\ell_1$ of Corollary \ref{c1} is a
direct application of the first part of the proposition. Now, let
us prove the irrationality of the constant $\ell_2$ appearing in
Corollary \ref{c2}. We must notice that this is not a direct
application of the first part of the proposition, because the set
of the prime numbers has asymptotic density $0 < 1 - \log 2$.

Let $\F$ denote the set of square-free numbers, that is the set of
all positive integers which are a product of pairwise distinct
prime numbers. It is known that $\F$ has asymptotic density
$\frac{6}{\pi^2} > 1 - \log 2$ (see e.g. \cite{hw}). So, it
follows by Erd\H{o}s' result that the number
$$\ell_3 := \sum_{n \in \mathbb{N}} \frac{1}{\lcm\{f \in \F ~|~ f \leq n\}}$$
is irrational.\\
But we remark that for all $n \in \mathbb{N}$, we have
$$\lcm\{f \in \F ~|~ f \leq n\} ~= \!\!\!\prod_{\begin{subarray}{c}
\scriptstyle{p ~\text{prime}} \\ \scriptstyle{p \leq n}
\end{subarray}} \!\!\!\!p ~=~ \lcm\{p ~\text{prime} ~|~ p \leq n\} ,$$
which shows that actually $\ell_3 = \ell_2$. Consequently $\ell_2$
is an irrational number. This completes the proof of the
proposition.\hfill $\blacksquare$~\vspace{2mm}

\noindent{\bf Remark.} The irrationality of both constants
$\ell_1$ and $\ell_2$ appearing in Corollaries \ref{c1} and
\ref{c2} respectively can be shown by a more elementary way than
that presented in Erd\H{o}s' paper for the general case.

\end{document}